\documentclass[12pt]{article}
\newcommand{\be}{\begin{equation}}
\newcommand{\bc}{\begin{center}}
\newcommand{\ec}{\end{center}}
\newcommand{\ee}{\end{equation}}
\newcommand{\bea}{\begin{eqnarray}}
\newcommand{\eea}{\end{eqnarray}}
\newcommand{\beaa}{\begin{eqnarray*}}
\newcommand{\eeaa}{\end{eqnarray*}}

\newcommand{\NN}{\mathbb{N}}

\newcommand{\RR}{\mathbb{R}}

\newcommand{\eps}{\varepsilon}
\newtheorem{theorem}{Theorem}
\newtheorem{lemma}{Lemma}

\newtheorem{remark}{Remark}

\newcommand{\pt}{\partial_t}

\newcommand{\norm}[1]{\left\|#1\right\|}

\newcommand{\ds}{\,\mathrm{d}s}

\newcommand{\dsp}{\displaystyle}

\begin{document}

%\begin{frontmatter}
%\pagestyle{headings}
\bc
{\bf \large On the long-time behavior of a continuous duopoly model with constant conjectural variation }

Isabella Torcicollo
%ON THE LONG-TIME BEHAVIOR OF A CONTINUOUS DUOPOLY MODEL WITH CONSTANT CONJECTURAL VARIATION

Istituto per le Applicazioni del Calcolo ``Mauro Picone'',
CNR, \\
Via P. Castellino 111, Naples, Italy. i.torcicollo@iac.cnr.it
\ec
{\bf Abstract}
The paper concerns a continuous model governed by a ODE system originated by a discrete duopoly model  with bounded rationality, based on constant conjectural variation. The ultimately boundedness of the solutions (existence of an absorbing set in the phase space) is shown and conditions guaranteeing the nonlinear, global, asymptotic stability of solutions have been found using the Liapunov direct method.


\if0
%
%
%% Place the running title of the paper with 40 letters or less in []
 %% and the full title of the paper in { }.
%\title
%[A continuous  reaction-diffusion duopoly model ] %Use the shortened version of the full title
      {A continuous  reaction-diffusion duopoly model with bounded rationality and constant conjectural variation:  action of  outputs self and cross diffusion
      \thanks{This paper has been performed under the auspices of G.N.F.M. of I.N.D.A.M.
 The author acknowledges the Project G.N.F.M. Giovani 2015 \lq \lq Dinamica di sistemi complessi con applicazioni in Fluidodinamica, Biologia ed Economia". } }

% Place all authors' names in [ ] shown as running head, Leave { } empty
% Please use `and' to connect the last two names if applicable
% Use FirstNameInitial.  MiddleNameInitial. LastName, or only last names of authors if there are too many authors

%The abstract of your paper

%\begin{AMS} 35B40, 35B35, 35B57;  34D23 \end{AMS}  \pagestyle{myheadings}  \thispagestyle{plain} \markboth{TEX PRODUCTION}{USING SIAM'S \LaTeX\ MACROS}
\fi
%%%%%%%%%%%%%%%%%%%%%%%%%%%%%%%%%%%%%%%%%%
%%%%%%%%%%%%%%%%%%%%%%%%%%%%%%%%%%%%%%%%%%%%%%%%%%%%%%%%%%%%%
\section{Introduction}

%%%%%%%%%%%%%%%%%%%%%%%%%%%%%%%%%%%%%%%%%%%%%%%%%%%%%%%%%%%%%%%%%%%%%%%%%%%%
The classic model of oligopoly was proposed by the French mathematician A. Cournot \cite{Cou} and dates back to 1838. The oligopoly market structure showing the action of only two companies  is called duopoly. Duopoly is an intermediate situation between monopoly and perfect competition, and analytically is a more complicated case, because an oligopolist must consider not only the behavior of the costumers, but also those of the competitors and their reactions. %This model has been studied in the literature from different points of view, see for instance [2?7] or [8].  
In duopoly game each duopolist believe that he can calculate the quantity he should produce in order to maximize his profits.
In the study of theoretic and realistic problems of duopoly, the Cournot model
% (Cournot was the first to propose the duopoly model in 1838 \cite{Cou})
 and successively the conjectural variation model 
proposed by Bowley in 1924 \cite{Bow} and later by Frish in 1933 \cite{Fri}
have been generally accepted and become the two important tools, convenctional game theory and conjectural variation model, to describe  the market behavior. A useful summary of the history and of the debate on the conjectural variations is provided by Giocoli in \cite{G}. Conjectural variation is referred to the beliefs that one firm has about the way its competitor may react if it varies its output. The firm forms a conjecture about the variation in the other firm's output that will accompany any change in its own output. In the classic Cournot model, it is assumed that each firm treats the output of the other firm as given when it chooses its output. In this case the conjectural variation is zero.\\
Many researchers  (see for instance  \cite{A}-\cite{BN}, \cite{E}, \cite{K}), %such as Kopel \cite{}, Bischi et al. \cite{}, Ahmed and Agiza [8], and Agiza and Elsadany [9], 
have paid a great attention to the dynamics of games. 
To make the theory more realistic, people have put forward various bounded rationality behaviors, combined game theory with dynamic system, assumed a delay in the term of bounded rationality, and formed this field of bounded rational dynamic game.
Different expectations such as naive expectations, adaptive expectations and bounded rationality have been proposed. Generalizations have been made for the Cournot-Kopel duopoly model \cite{T},\cite{RT1}, by introducing self and cross diffusion terms in the equations \footnote{The reaction-diffusion systems of PDEs are the best candidates for investigating the diffusion processes in spatial domains (see for instance \cite{CDT}-\cite{CDDT}, \cite{DA}, \cite{10}, \cite{R}-\cite{Ri}, \cite{RT}- \cite{TV3}).}, when the outputs are in the market in large territories.\\
%The theoretical development of complex duopoly dynamics has been recently surveyed in [10, 11].
 %Oligopoly theory is one of the oldest branches of mathematical economics  and dates back to 1838 when its basic model was proposed by Cournot \cite{Cou}.
%The first well-known model which gives a mathematical description of competition in a duopoly market dates back to the French economist Cournot (1838) \cite{Cou}.\\
%In fact, the first economist to introduce the idea that, to the duopolists' reaction functions had to be given a conjectural interpretation,  is Bowley  in 1924 \cite{Bow}. The concept of conjectural variation, however, must wait until the work of Frish in 1933 \cite{Fri} to obtain a final theoretical arrangement. A useful summary of the history and of the debate on the conjectural variations is provided by Giocoli in \cite{G}.
%To make the theory more realistic, many authors have combined game theory with dynamic systems theory and formed the  field of bounded rational dynamic game (the general formula of duopoly game with a form of bounded rationality has been given in \cite{BN}). Some authors have assumed a delay in the term of bounded rationality \cite{ZZ}, \cite{E}, while others   have proposed many types of expectations (a  player, in order to adjust the own output, can choose the strategy rule among many available techniques). Naive or adaptive expectations,  boundedly  rational strategies are just  a small example and, consequently, many models have been constructed (see for instance \cite{A}-\cite{BK},  \cite{EAE}-\cite{K} and references therein). \\
Given a model, economists need to make predictions on the asymptotic behavior of the system i.e. how the model behaves in the future. %n order to study the longtime behavior, the mathematical models play a central role. In fact, 
Although the mathematical models represent only an approximation of the problem (they consider only some variables that are involved in the phenomenon), they allow to obtain an estimate for the market development. \\
Recently, in \cite{YY}, the  discrete duopoly model (\ref{br7}) with bounded rationality based on constant conjectural variation has been considered and the local stability of critical points, numerical  simulations,   Liapunov exponents and fractal dimensions of strange attractors are investigated. 
% In addition, symmetry model and Bertrand model have been investigated. 
 In the present paper, we generalize (\ref{br7})  via a continuous ODEs system aimed to investigate on the longtime behavior of the solutions. Precisely, our aim is to  
 \begin{itemize}
 \item[i)] show that, in the phase space, there exists an absorbing set of the solutions;
 \item[ii)] find conditions assuring the nonlinear global asymptotic stability of solutions.
 \end{itemize}
  %, disregarding the diffusion in the territory in which the outputs are in the market \footnote{Misregarding the diffusion processes, can not  be allowed when the territories in which the outputs are in the market are sufficiently large.}.  
% We show that, the long-time behavior of the ODEs model is essentially the same of the discrete model.
% Successively, we take into account the outputs self and cross diffusion.  
 %Since the reaction-diffusion systems of PDEs are the best candidates for investigating the diffusion processes in spatial domains   \{see  \cite{Cap} - \cite{CGH}, \cite{DA} - \cite{10}, \cite{R} - \cite{Sor}, \cite{RT} - \cite{RT3}, \cite{T1} - \cite{TV2}
%\cite{D}, \cite{Ri}-
%and references therein\},  when the outputs are in the market in large territories, we generalize (\ref{br9}) via the reaction-diffusion PDEs model (\ref{dif}) with self and cross diffusion \footnote{We remark that a similar generalization has been done in \cite{T},\cite{RT1} for the Cournot-Kopel duopoly  model, with adaptive expectations.}. We show, (as expected),  that the action of the outputs self and cross diffusion can be very relevant.  Precisely, a steady state stable (unstable)  according to the discrete and ODEs model, can become unstable (stable) for large sets of values of the diffusion coefficients $\gamma_{ij}.$\\
The paper is organized as follows. It begins by recalling the discrete model (\ref{br7}) proposed in \cite{YY} (Section 2). Successively, the continuous counterpart of (\ref{br7}) is introduced and the invariance of the first orthant is shown (Section 3).  Section 4 is devoted to the existence of an absorbing set and uniqueness. Critical points are determined in Section 5 while linear and nonlinear global asymptotic stability are studied respectively in Section 6 and 7. 
%All the subsequent Sections are devoted to (\ref{dif}). Precisely, in Section 8, (\ref{dif})  is introduced and the functional space in which is embedded is precised, while in Section 9,  the linear stability of the constant steady states  is considered with respect to perturbations deviluppable in Fourier series of a complete sequence of eigenfunctions of (\ref{spectral}) with positive eigenvalues. In Section 10 the role of the self and cross diffusion on the linear stability is put in evidence, while, in Section 11, an interpretation of stability-instability influence on  spreadings of outputs is given. In Section 12, the nonlinear local and global stability of $E_3$ is investigated. Section 13 is devoted to some final remarks and comments. 
The paper ends with an Appendix (Section 8) in which the origin of the function (\ref{functional}) is recalled. %(Subsection 8.1) 
%and sketches  of proof of Lemma 1  (Subsection 14.2) and of Theorem 12 (Subsection 14.3) are given. 
%%%%%
%%%%%%%%%%%%%%%%%%%%%%%%%%%%%%%%%%%%%%%%%%%%%%%%%%%%%%%%%%%%%
\section{The discrete model}
%%%%%%%%%%%%%%%%%%%%%%%%%%%%%%%%%%%%%%%%%%%%%%%%%%%%%%%%%%%%%%%%%%%%%%%%%%%%
%%%%%%%%%%%%
%%%%%%%%%%%%%%%%%%
%%%%%%%%%%%%%%
In the duopoly model, proposed in \cite{YY}, it is assumed that two duopolists $X$ and $Y,$ in the hypothesis of constant conjectural variation, produce similar products for quantity competition. The firms do not have complete information or complete rationality behavior, but they know marginal profit function of constant conjectural variation. It is a model  which utilize marginal profit function of static constant conjectural variation model to make output decision.
The marginal profit function of conjectural variation of firms $X$ and $Y$ has the following form
\be \left\{ 
\begin{array}{l} 
 \label{1}
\dsp \Pi_{\small x}(x,y)= \theta_1-\gamma y-L_1x\\
\dsp \Pi_{\small y}(x,y)= \theta_2-\gamma x-L_2 y
\end{array}\right.
\ee
with $\theta_i=\alpha_i-c_i>0$ (i=1,2) where $c_i>0$ is the marginal cost function, while $\alpha_i,$  $L_i$ and $\gamma$ are positive model parameters. Denoted by $x_t$ and $y_t$ the output of firm $X$ and $Y$ respectively, the firms adjust the output of the next period $t+1$ considering marginal profit function and output of the current period $t$.  When marginal profit  function is positive they have an action of increasing output; when marginal profit function is negative they have an action of decreasing output. One obtains the following discrete system for the two firms
  \be  \left\{ \begin{array}{ll}
x_{t+1}=x_t + a x_t (\theta_1- \gamma y_t - L_1 x_t)  \label{br7}
\\
y_{t+1}=y_t + \nu y_t (\theta_2- \gamma x_t - L_2 y_t)   \end{array}\right.\ee
where $a>0$ and $\nu>0$ represent the speeds of output adjustment of two firms.
In  \cite{YY},  the local stability of critical points, numerical  simulations,  bifurcation diagram,  Liapunov exponents and fractal dimensions of strange attractors are investigated.  In addition, symmetry model and Bertrand model have been investigated. 
%
%
%%%%%%%%%%%%%%%%%%%%%%%%%%%%%%%%%%%%%%%%%%%%

%%%%%%%%%%
%%%%%%%%%
%%%%%%%%%
%%%%%%%%%%%%%%%%%%%%%%%%%%%%%%%%%%%%%%%%%%%%%%%%%%%%%%%%%%%%%%%%%%%
\section{Continuous ODEs model}
%%%%%%%%%%%%%%%%%%%%%%%%%%%%%%%%%%%%%%%%%%%%%%%%%%%%%%%%%%%%%%%%%%%
Assuming continuous time scales, denoting by $u$ and $v$ the outputs of the two firms $X$ and $Y$, respectively, from (\ref{br7}) the nonlinear continuous system for the evolution of $u$ and $v$ is immediately 
obtained \be \left\{ \begin{array}{l}\dsp \frac{du}{dt}=a u(\theta_1- \gamma v - L_1 u)  \label{br9}\\
\\
\dsp \frac{dv}{dt}=\nu v(\theta_2- \gamma u - L_2 v) 
 \end{array} \right.
\ee
where
%in general,  $a$, $\nu$, $\gamma$, $\theta_i$ and $L_i$ with  $(i=1,2)$ are positive model parameters,  
%$\Omega \subset \RR^3,$  is a bounded domain (territory),
%\footnote{It has been considered $\Omega \subset \RR^3$ also for taking into account territories containing mountains.},
%in which, the outputs are in the market, 
\[ \phi: t \in  {R}^+ \rightarrow \phi (t) \in
R, \qquad \phi\in \{u,v\}.\]
%$$v:({\bf x}, t)\in \Omega \times{\RR}^+ \rightarrow v({\bf x}, t) \in \RR,$$ 
 To (\ref{br9}) we
append the initial data 
\be  \label{inn}
 u( 0)=u_0,   ~~~~~~~~
v(0)=v_0 \ee
with $u_0,v_0$  assigned positive constants.\\
Further we remark that the first orthant is invariant. In fact, integrating (\ref{br9}), one obtains 
\be \left\{\begin{array}{l}
\dsp u=u_0\exp{\dsp\int_0^t a(\theta_1-\gamma v -L_1 u)d\tau}\\
\dsp v=v_0 \exp{\dsp\int_0^t \nu(\theta_2-\gamma u -L_2 v)d\tau}
\end{array} \right.\ee
and hence $\{u_0>0, v_0>0\}$ imply $\{u(t)>0, v(t)>0  \, \, \forall t\geq 0\}$  i.e. the first orthant is invariant.
\section{Ultimately boundedness (absorbing sets) and uniqueness}
As it is well-known a set ${\mathcal A}$ of the phase space $(u,v)$, is an absorbing set if, denoting by $d[(u,v),{\mathcal A}]$ the distance between $[u(t),v(t)]$ and ${\mathcal A}$ at time $t,$ i.e. 
\be
d(t)= \inf_{\mathcal A}(\vert u- \bar U\vert^2 + \vert v-\bar V\vert^2),  \qquad (\bar U, \bar V) \in {\mathcal A},
\ee
it follows that 
\begin{itemize}
\item[i)] ${\mathcal A}$ is a global attractor, i.e. 
\be
\lim_{t\rightarrow \infty} d[(u,v), {\mathcal A}] =0,
\ee
with $(u_0,v_0) \in {\mathcal B}$, for any open set ${\mathcal B}\supset {\mathcal A};$
\item[ii)] ${\mathcal A}$ is positively invariant, i.e.
\be
(u_0,v_0) \in {\mathcal A} \Rightarrow (u(t), v(t)) \in {\mathcal A}, \quad \forall t >0.
\ee
\end{itemize}
\begin{theorem}
Any set containing the rectangle 
\be
S=\left\{ (u,v) \in R^2_+ : 0<u\leq\dsp\frac{\theta_1}{L_1},\,\,  0<v\leq\dsp\frac{\theta_2}{L_2}  \right\},
\ee
of the phase space is an absorbing set for (\ref{br9}).
\end{theorem}
{\em Proof.}  By virtue of  $(\ref{br9})_1$ it follows that
\be 
\dsp \frac{du}{dt}=a u(\theta_1- \gamma v - L_1 u) \leq a u(\theta_1-  L_1 u) 
\ee
from which, by setting $w=\dsp \frac 1u$, we obtain
\be \label{eqd}
\dsp \frac{dw}{dt}\geq -a \theta_1w +a L_1 ,
\ee
and hence, integrating, it follows that 
\be \label{uf}
w\ge w_0  e^{-a\theta_1 t} + \frac{L_1}{\theta_1}(1-e^{-a\theta_1 t}).
\ee
Since (\ref{uf}) implies 
\be
\dsp\lim_{t\rightarrow \infty} w\ge \frac{L_1}{\theta_1},
\ee 
one immediately obtains that,  $\forall \varepsilon>0$, there exists $T_1(\varepsilon)>0$ such that 
\be \label{ut}
 u(t) \le \frac{\theta_1}{L_1}+\varepsilon ,  \qquad \forall t \geq T_1.
\ee
Following the same procedure for (\ref{br9})$_2,$ one obtains the existence of a  constant $T_2(\varepsilon)>0$   such that 
\be  \label{vt}
 v(t) < \frac{\theta_2}{L_2}+\varepsilon ,  \qquad \forall t \geq T_2.
\ee
For $\varepsilon \rightarrow 0$ and $T=\max\{ T_1,T_2\},$ from (\ref{ut})-(\ref{vt}),  it follows that
\be \label{uvt}
 u(t) \leq \frac{\theta_1}{L_1},  \qquad 
 v(t) \leq \frac{\theta_2}{L_2},  \qquad \forall t \geq T,
\ee
i.e. S is an attractor.\\
In order to prove that S is positively invariant, let us assume that
\be \label{uvz}
 0<u_0 < \frac{\theta_1}{L_1},  \qquad 
 0<v_0 < \frac{\theta_2}{L_2},
\ee
then, from  (\ref{uf}) and  (\ref{uvz})$_1$ it follows that 
\be \label{ufu}
u(t)\leq
% \dsp \frac{u_0 \theta_1}{\theta_1 e^{-a\theta_1 t} + u_0L_1(1-e^{-a\theta_1 t})} \leq 
 \dsp \frac{u_0 \theta_1}{u_0 L_1 e^{-a\theta_1 t} + u_0L_1(1-e^{-a\theta_1 t})} =\dsp \frac{ \theta_1}{L_1} \,\, \, \, \,\, \,  \forall t\geq T.
\ee
Following the same steps for $v(t)$, one obtains 
\be \label{vfv}
v(t)\leq \dsp \frac{ \theta_2}{L_2} \,\, \, \, \,\, \,  \forall t\geq T,
\ee
i.e. S is invariant.\\
\begin{remark} Since any solution enters in S in a finite time and remains there, without loss of generality, we can confine ourselves to study the dynamical behavior of the system in S assuming $(u_0,v_0)\in S$.
\end{remark}
%\end{proof}\\
\noindent The following uniqueness theorem holds
\begin{theorem}
The system (\ref{br9}) - (\ref{inn})  admits a unique solution.
\end{theorem}
{\em Proof.} Let $(u_1,v_1)$  and $ (u_2, v_2)$ be two solutions of (\ref{br9})-(\ref{inn}).
On setting
\be
\tilde u=u_1-u_2, \,\,\, \,\,\, \tilde v=v_1-v_2
\ee
it follows that 
\be\label{UV}
\left\{\begin{array}{l}
\dsp\dsp \frac{d\tilde u}{dt} =a \theta_1 \tilde u - a L_1( u_1+ u_2)\tilde u -a \gamma (u_1 v_1 - u_2 v_2) 
\\
\\
\dsp\dsp \frac{d\tilde v}{dt} =\nu \theta_2 \tilde v - \nu L_2( v_1+ v_2)\tilde v - \nu  \gamma (u_1 v_1 - u_2 v_2) 
  \end{array}\right.\ee
with 
\be\label{intilde}
%\left\{\begin{array}{ll}
 \tilde u (0)=0,\,\,\,\,\, \,\, 
  \tilde v (0)=0.
%\end{array}\right.
\ee
 From (\ref{UV}) it follows that
\be\label{calEtilde}
\begin{array}{l}
\dsp \frac12 \frac{d}{dt}({\tilde u}^2+{\tilde v}^2)= a[\theta_1 -L_1(u_1+u_2)]\tilde u^2 + \nu [\theta_2-L_2(v_1+v_2)]\tilde v^2 \\
\\
-a\gamma ( u_1  v_1- u_2 v_2)   \tilde u -\nu \gamma ( u_1 v_1 - u_2 v_2) \tilde v .
  \end{array}\ee
By virtue of the boundedness of solutions and $\vert \tilde u\vert +\vert \tilde v\vert < \sqrt 2 (\tilde u^2 +\tilde v^2)^\frac12$   (\ref{calEtilde}) leads to
\be\label{fintilde}
\begin{array}{l}
\dsp  \frac{d}{dt}({\tilde u}^2+{\tilde v}^2)\leq \kappa(\tilde u^2 + \tilde v^2 )^{\frac12} [1+(\tilde u^2 + \tilde v^2 )^{\frac12} ],
  \end{array}\ee
with $\kappa$ positive constant. Integrating (\ref{fintilde}) it  immediately follows 
\be
\dsp{\tilde u}^2+{\tilde v}^2\leq\frac{ ({\tilde u_0}^2+{\tilde v_0}^2)e^{\kappa t}}{[1+({\tilde u_0}^2+{\tilde v_0}^2)^\frac12(1-e^{\kappa t/2})]^2}   \qquad \forall t\geq 0.
\ee
and in view of (\ref{intilde}), one obtains
 \[
{\tilde u}^2+{\tilde v}^2\equiv 0.
\]
%\end{proof}
%%%
%%%
%%%
\section{Critical points of (\ref{br9})}
%
%%%%%%%%%%%%%%%%%%%%%%%%%%%%%%%%%%%%%%%%%%%%%%%%%%%%%%%%%%%%%%%%%%%
The critical points of (\ref{br9}) are the roots of the algebraic system
 \be
\left\{ \begin{array}{l}
 a u (\theta_1- \gamma v - L_1 u)=0   \label{alg1}
\\
 \nu v (\theta_2- \gamma u - L_2 v)=0 . 
\end{array} \right.
\ee
%We denote by  $(\bar u,\bar v)$ the generic equilibrium point of (\ref{br9}). 
Besides  the trivial equilibrium $E_0=(0,0) $,  system (\ref{br9})  admits the following nontrivial
critical points (boundary equilibrium points)
\be \label{bound}
E_1=(0, \dsp \frac{\theta_2}{L_2}), \qquad E_2=(\dsp \frac{\theta_1}{L_1},0 ) %\qquad  {\rm boundary \, equilibrium \,  points}
\ee
and  the only equilibrium of constant conjectural variation 
\be  \label{conj}
E_3=\left(\dsp \frac{\theta_1 L_2- \theta_2 \gamma}{L_1L_2-\gamma^2}, \, \,   \dsp \frac{\theta_2 L_1- \theta_1 \gamma}{L_1L_2-\gamma^2} \right),
\ee
where the parameters satisfy 
\be \left\{
\begin{array}{l}  \theta_1 L_2- \theta_2 \gamma>0
\\
\theta_2 L_1- \theta_1 \gamma>0
\\
L_1L_2-\gamma^2>0 .  \label{posE}
 \end{array} \right.\ee
\begin{remark}
The property of the equilibrium  of conjectural variation is that:
\begin{itemize}
\item[i)] when firm 1  chooses equilibrium $\bar u,$  then $\bar v$  is the optimal selection of firm 2;
\item[ii)]  when firm 2 chooses $\bar v$, then $\bar u$ is the optimal selection of firm 1;
\item[iii)] the realization condition of equilibrium of conjectural variation is that each firm has complete information of market and the ability of complexity rational behavior.
\end{itemize} 
\end{remark}
\section{Linear stability}
Denoting by $(\bar u, \bar v)$ a constant critical point of (\ref{br9}), we set
 \be U=u-\bar u,~~~~~~ V=v-\bar v.\label{pertu}\ee  
Then, the  equations  governing $(U,V)$ are given by 
\be
\left\{\begin{array}{ll} 
\dsp\frac{ dU}{dt} = 
a (\theta_1- \bar v \gamma  -2 \bar u L_1 )U - a \bar u \gamma V - a \gamma UV- a L_1 U^2\\
\\
\dsp\frac{ dV}{dt} = -\nu \gamma \bar vU + \nu (\theta_2- \gamma \bar u - 2 \bar v L_2 )V- \nu \gamma UV- \nu L_2 V^2
.\label{din}
\end{array}\right.\ee
Setting 
  \be\left\{\begin{array}{ll}
f(U,V)=- a \gamma UV- a L_1 U^2     ~~ ~~~   & g(U,V)=- \nu \gamma UV- \nu L_2 V^2\\
a_{11}=a (\theta_1- \bar v \gamma  -2 \bar u L_1 )    ~~   &a_{12}= - a \bar u \gamma  \\
 a_{21}= -\nu \gamma \bar v   ~~~~    & a_{22}= \nu (\theta_2- \gamma \bar u - 2 \bar v L_2 )    \label{position}
\end{array}\right.
 \ee
 (\ref{din}) becomes
\be\left\{\begin{array}{ll}\dsp\frac{ dU}{dt} =
a_{11} U +a_{12}  V +f(U,V)\\
\\
\dsp\frac{ dV}{dt} = a_{21}  U + a_{22} V + g(U,V) ,\label{perinc}
\end{array}\right.
\ee
under the  initial data
\be\label{init}
%\left\{\begin{array}{ll}
 U(0)=u_0-\bar u, \,\,\,\,\,\,\,
  V(0)=v_0-\bar v.
%\end{array}\right.
\ee
Linearizing, one obtains
\be \label{linearL}
\dsp \frac{d}{dt} \left(\begin{array}{l}U\\
V
\end{array}
  \right) = {\cal L}  \left(\begin{array}{l}U\\
V
\end{array}
  \right) 
\ee
with ${\mathcal L}$  given by 
\be
{\cal L}=  \left(\begin{array}{ll}a_{11} \,\,\,&a_{12}\\
a_{21} & a_{22}
\end{array}
  \right).
\ee
Denoting by $\lambda_i \, (i=1,2)$ the eigenvalues of ${\mathcal L},$ it is well-known that 
\be \left\{
\begin{array}{l} I_0=a_{11}+a_{22}=\lambda_1+\lambda_2\\
A_0=a_{11}a_{22}-a_{12}a_{21}=\lambda_1\lambda_2
\label{condii}
\end{array}\right.
\ee 
 and hence, {\it if and only if,
 \be \left\{
\begin{array}{l} I_0(\bar u, \bar v)<0\\
A_0(\bar u,\bar v)>0
\label{aut}
\end{array}\right.
\ee 
the critical point $(\bar u, \bar v)$ is linearly stable.}
\begin{theorem}
By virtue of (\ref{aut}), it follows that
% For the linear system (\ref{linearL}), one obtains that 
\begin{itemize}
\item[i)] $E_3$ is linearly stable;
\item[ii)] $E_0, E_1, E_2$  are unstable;
\item[iii)] $E_0$ is a unstable node; $E_1,E_2$  are  saddle point; $ E_3$ is a stable node;
\item[iv)] the spectral equation is given by  $\lambda^2-I_0 \lambda+ A_0=0;$  
\item[v)]  $(I_0^2-4A_0)_{(\bar u, \bar v)=E_i}>0, \,\,\forall i=\{0,1,2,3\}$ i.e. the eigenvalues are real. 
\end{itemize}
\end{theorem}
{\em Proof.} By virtue of (\ref{aut}), one obtains 
%In order to prove {\it  i) ii)}, let us apply (\ref{condii})-(\ref{aut}) to $E_i (i=0,1,2,3)$. In view of (\ref{posE}), it easily %follows that 
\be \left\{\begin{array}{ll}  \label{E33}
I_0(E_3)= - \dsp \frac{aL_1(\theta_1 L_2- \theta_2 \gamma)}{L_1L_2-\gamma^2}- \dsp \frac{\nu L_2(\theta_2 L_1- \theta_1 \gamma)}{L_1L_2-\gamma^2} <0
\\
A_0(E_3)= \dsp \frac{a\nu(\theta_1 L_2- \theta_2 \gamma)(\theta_2 L_1- \theta_1 \gamma)}{L_1L_2-\gamma^2}>0,
\end{array}  \right.
\ee
\bigskip
\be \left\{\begin{array}{ll} \label{E00}
I_0(E_0)=a\theta_1+\nu\theta_2>0
\\
A_0(E_0)= a\theta_1 \theta_2\nu>0,
\end{array}  \right.
\ee
\be \left\{\begin{array}{ll} \label{E11}
I_0(E_1)=  \dsp\frac{ a(\theta_1 L_2- \theta_2 \gamma)}{L_2}- \nu\theta_2 
\\
A_0(E_1)= \dsp\frac{ a(\theta_1 L_2- \theta_2 \gamma)}{L_2}(- \nu\theta_2 )<0,
\end{array}  \right.
\ee
\be \left\{\begin{array}{ll}  \label{E22}
I_0(E_2)=  \dsp\frac{ \nu(\theta_2 L_1- \theta_1 \gamma)}{L_1}- a\theta_1
\\
A_0(E_2)= \dsp\frac{ \nu(\theta_2 L_1- \theta_1 \gamma)}{L_1}(-a\theta_ 1)<0,
\end{array}  \right.
\ee
and hence,  {\it  i)- ii)} immediately follow.\\
In order to prove {\it iii), iv), v)}, from
\be
det \left(\begin{array}{ll}a_{11} -\lambda \,\,\,&a_{12}\\
a_{21}& a_{22}-\lambda 
\end{array}
  \right) =
 \lambda^2-(a_{11}+a_{22})\lambda+(a_{11}a_{22}-a_{12}a_{21})=0
\ee
it follows that, the spectral equation is given by $\lambda^2-I_0\lambda+A_0=0.$\\
In addition, since $$I_0^2-4A_0= (a_{11}-a_{22})^2+4a_{12}a_{21}>0\,\,\,\,\,\,\,\forall  E_i \,\, (i=0,1,2,3)$$ 
%_{(\bar u, \bar v)=E_i}>0,  \forall i=\{0,1,2,3\}$$  i.e. the eigenvalues are real. 
the eigenvalues are real and {\it v)} holds. \\
Further, it immediately  follows that  \\
$E_0$ is an unstable node since  $\left\{\lambda_1=  a\theta_1 >0, \,\,\,
\lambda_2= \nu\theta_2>0  \right\}$;\\
\\
$E_1$ is a saddle point, since $\left\{\lambda_1= \dsp\frac{ a(\theta_1 L_2- \theta_2 \gamma)}{L_2}>0, \,\,\,
\lambda_2=- \nu\theta_2 <0 \right\};$\\
\\
$E_2$ is a saddle point, since $\left\{\lambda_1= - a\theta_1 <0,\,\,\,
\lambda_2=\dsp\frac{ \nu(\theta_2 L_1- \theta_1 \gamma)}{L_1}>0 \right\};$\\
$E_3$ is a stable node, since  $\left\{ \begin{array}{l} \lambda_1=-\vert I_0\vert  - \dsp \sqrt{I_0^2 -4A_0}  <0, \\
\lambda_2=-\vert I_0\vert  + \dsp \sqrt{I_0^2 -4A_0}  <0 \end{array} \right\}.$\\
%\end{proof}
%
%
%
\section{Nonlinear stability of $E_3$}
The nonlinear stability analysis is based on the Rionero's function  \cite{Dir},  \cite{rion}-\cite{ternary},
\be \label{functional}
{\mathcal V} = \dsp\frac12 [A_0( U ^2 +  V ^2) + {(a_{11}V- a_{21}U)}^2+ {( a_{12} V- a_{22} U)}^2]
\ee
whose time derivative along the solutions of  (\ref{perinc}) is given by 
\be  \label{dotv}
\dsp \dot{\mathcal V}= -A_0\vert I_0\vert ( U^2+ V^2) + \Psi
\ee
with  $A_0$, $I_0$ given by  (\ref{condii}),
\be  \label{Psi}
\Psi =( \alpha_1 U-\alpha_3 V) f( U, V) + ( \alpha_2 V-\alpha_3 U)  g( U, V) 
\ee
and 
\be  \begin{array}{ll}
%\left\{ \begin{array}{ll}
\alpha_1= A_0+ a_{21}^2+a_{22}^2>0, 
%\, \,\,\,\,\, 
&\alpha_2=A_0+a_{11}^2+ a_{12}^2>0, \\
\alpha_3=  a_{11}a_{21}+a_{12}a_{22}>0 .\,&
\end{array}
\ee
\begin{theorem} 
The equilibrium of conjectural variation $E_3$ is  nonlinearly asymptotically exponentially  stable, according to 
\be  \label{exp}
%{\cal V} \leq {\cal V}(0)  \exp[-(h_1- h_2{\cal V}^{\frac12}(0))t]  \qquad  \forall t>0
{\mathcal V} \leq {\mathcal V}(0) e^{-(1-\eta) h_1 t}   \qquad  \forall t>0
\ee
 with $\eta\in ]0,1[$ and  $h_1$   positive constant,
% $(U,V)$ be a solution to  (\ref{din}) under conditions  (\ref{init})- (\ref{robin})
 under the local condition 
\be
U_0^2+V_0^2 \leq \dsp \frac{(A_0 \vert I_0\vert \delta_1 )^2}{2M^2\delta_2^2}
\ee
with $\delta_1,\delta_2$ positive constants 
and globally nonlinearly asymptotically exponentially stable under the condition 
\be \label{ass}
\dsp\frac{\theta_1^2}{L_1^2}+ \frac{\theta_2^2}{L_2^2} 
 %\dsp \left(\frac{(\theta_2L_1-\theta_1\gamma)\gamma}{L_1(L_1L_2-\gamma^2)}\right)^2+ \dsp %\left(\frac{(\theta_1L_2-\theta_2\gamma)\gamma}{L_2(L_1L_2-\gamma^2)}\right)^2
\leq \dsp \frac{(A_0 \vert I_0\vert \delta_1 )^2}{2M^2\delta_2^2}.
\ee 
%\be
%(\vert U_0 \vert^2 + \vert V_0 \vert^2 ) ^\frac12 < \frac{\bar \lambda}{\sqrt 2 M},
%{\cal V}^\frac12(0) <\dsp \frac{h_1}{h_2}
%\frac{A_0\vert I_0\vert \delta_1}{ \sqrt 2 M \delta_2}
%\label{picdat1}
%\ee 
% \be  \label{M}  M= \max (a\gamma, \nu \gamma, aL_1, \nu L_2)  \ee
\end{theorem}
{\it Proof.}   From  (\ref{Psi}) immediately one obtains 
\be 
\begin{array}{l} \label{nlt1}
\dsp  [\alpha_3 ( a L_1 +\nu \gamma)-\alpha_1  a \gamma] U^2V- \alpha_1\ a L_1U^3 \\
\\
+\dsp  [\alpha_3( \nu L_2+a\gamma) -  \alpha_2 \nu \gamma]UV^2 
- \alpha_2 \nu L_2 V^3\\
\\
\dsp \leq  M (U^2+V^2)(\vert U\vert +\vert V\vert) \leq 
\sqrt 2 M  ({U}^2+{V}^2)^\frac32 ,
  \end{array}
\ee
where 
\[
\begin{array}{ll}
M= \max (m_1,m_2,m_3, m_4) \,\,\, & m_1= \vert \alpha_3( a L_1 +\nu \gamma)-\alpha_1 a \gamma\vert \\
m_2=\vert \alpha_3( \nu L_2+a\gamma) -  \alpha_2\nu \gamma\vert \,\,\,& m_3=\alpha_1 aL_1, \,\,\,\,\,\,\, m_4= \alpha_2\nu L_2. 
\end{array} 
\] 
Setting  $E=( U^2+  V^2)$ and on taking into account that 
%${\cal V}$ and $E$  are $L^2$-equivalent norms, i.e. 
 there exist two positive constants $\delta_1,\, \delta_2$ such that 
\be
\delta_1E \leq  {\mathcal V} \leq \delta_2 E
\ee
with $\dsp \delta_1=\frac12 A_0$ and $\dsp\delta_2= \frac12 A_0+ a_{11}^2+ a_{21}^2+a_{12}^2 + a_{22}^2,$
in view of   (\ref{nlt1}), the  (\ref{dotv})  leads to 
\be
\dot{\mathcal V} \leq -\dsp \frac{A_0\vert I_0\vert}{\delta_2}{\mathcal V} + \frac{\sqrt 2 M}{\delta_1^{3/2}}  {\mathcal V}^{\frac32}
\ee
and hence
\be
\dot{\mathcal V} \leq -(h_1- h_2{\mathcal V}^{\frac12}){\mathcal V},
\ee
where $h_1=\dsp \frac{A_0\vert I_0\vert }{\delta_2}, \,\,\,h_2=\dsp \frac{\sqrt 2 M}{ \delta_1^{3/2}}$. Choosing
\be
h_2 {\mathcal V}^\frac12(0) < h_1, \,\,\,{\rm i.e. }\,\,\,\,\,\,  h_2 {\mathcal V}^\frac12(0) =\eta  h_1, \,\,\,\,\,\,\,\, \eta\in ]0,1[
\ee 
%that is  $h_2 {\cal V}^\frac12(0) =\eta  h_1$
 by recursive argument,  (\ref{exp}) immediately follows.\\
On the other hand, in view of Theorem 1, one has 
\be
U_0^2+V_0^2 \leq  \dsp \frac{\theta_1^2}{L_1^2}+ \dsp \frac{\theta_2^2}{L_2^2}
\ee
and hence (\ref{ass}) immediately follows.
%\end{proof}
\\
\begin{remark}
In view of  i) and ii) of Theorem 3, it follows that, in the passage from the discrete model to the continuous ODEs one (with  constant conjectural variation), the constant steady states do not vary. Their stability-instability don't vary, except the fact that in the continuous model, $E_3$ is reached exponentially when $t \rightarrow \infty$, from any initial state. In fact, we remark that, since in the case $E\in \{E_0,E_1,E_2\}$, generally $I_0(E)A_0(E)<0$ does not continues to hold, therefore does not appear a stabilizing effect of the nonlinearities.
\end{remark}

\section{Appendix} 
%\subsection{ The basic peculiar Liapunov function}
We define \emph{peculiar} the Liapunov function able to give for the nonlinear stability-instability, exactly
the same conditions of the linear stability. We, following \cite{ternary}, here recall the construction of a such Liapunov function for a binary system of O.D.Es, since
is this function  that appears in Section 7.\\
Let us consider the stability of the zero solution of the system
\be\label{76}
\left\{\begin{array}{l}
\dsp\frac{dx}{dt}=ax+by+f(x,y),\\\\
\dsp\frac{dy}{dt}=cx+dy+g(x,y),
       \end{array}
 \right.
\ee
($f,g$ being nonlinear and such that $f(0,0)=g(0,0)=0$) and introduce the function
\be\label{77}
W=\dsp\frac{1}{2} I[A(x^2+y^2)+(ay-cx)^2+(by-dx)^2],
\ee
with
\be\label{78}
I=a+d=\lambda_1+\lambda_2,\,\,A=ad-bc=\lambda_1\cdot\lambda_2,
\ee
$\lambda_1,\lambda_2$ eigenvalues of
$\left(\begin{array}{ll}
a&b\\
c&d
       \end{array}
\right).$ Since 
\be\label{79}
\left\{\begin{array}{ll}
x\dot{x}=ax^2+bxy+xf,&y\dot{y}=cxy+dy^2+yg,\\\\
y\dot{x}=axy+by^2+yf,&x\dot{y}=cx^2+dxy+xg,
       \end{array}
\right.\ee
by straightforward calculations it follows that
\be\label{80}
\dsp\frac{dW}{dt}= I A(x^2+y^2)+\Psi,
\ee
with
\be\label{81}
\left\{\begin{array}{l}
\Psi=\eps I[(\alpha_1x-\alpha_3 y)f+(\alpha_2 y-\alpha_3 x)g],\\\\
\alpha_1=A+c^2+d^2,\,\,\alpha_2=A+a^2+b^2,\,\,\alpha_3=ac+bd.
       \end{array}
\right.\ee

\section*{Acknowledgments}
\noindent This paper has been performed under the auspices of G.N.F.M. of INdAM.
 The author acknowledges the Project G.N.F.M. Giovani 2015 \lq \lq Dinamica di sistemi complessi infinito dimensionali con applicazioni in Fluidodinamica, Biologia ed Economia".
 The author thanks gratefully Prof. Salvatore Rionero for his helpful suggestions and for his continuous encouragement and support.
%%

% \medskip
% The data information below will be filled by AIMS editorial staff
%Received xxxx 20xx; revised xxxx 20xx.
%\medskip

\begin{thebibliography}{0}
\bibitem{A} \newblock{H. N.  Agiza,} \newblock{On the analysis of stability, bifurcation, chaos, and
 chaos control of Kopel map,}  \emph{ Chaos Solitons Fractals},  {\bf 10} (1999) 1909--1916.
\bibitem{AEK} \newblock{H.N. Agiza, A.A. Elsadany and M. Kopel, } 
\newblock{Nonlinear dynamics in the Cournot duopoly game with heterogeneous
players,} \emph{ Physica A}  {\bf 320} (2003) 512--524.
\bibitem{BK}  \newblock{G.I., Bischi, L. Gardini, and M. Kopel,}  \newblock{Noninvertible maps and Complex basin boundaries in dynamic economic models with coexisting attractors,} \emph{ Chaos and Complexity Letters }, {\bf 2} (1)  (2006) 43--74.
\bibitem{BN} \newblock{G. I.,  Bischi and  A. Naimzada,}  \newblock{Global analysis of a dynamic duopoly game with bounded rationality,}  \emph{Advanced in Dynamic Games and Application,} {\bf 5}  Chapter 20, Birkhouser, 1999. 
\bibitem{Bow} \newblock{A.L. Bowley, }  \emph{The Mathematical Groundwork of Economics}. \newblock{Oxgord University Press. 1924.}
%\bibitem{Cap} \newblock{F. Capone, } \newblock{On the dynamics of predator-prey models with the Beddington-De Angelis functional response, under Robin boundary conditions,} \emph{Ricerche di Matematica} {\bf 57}  (2008), 137--157.
\bibitem{CDT} \newblock{F. Capone, R. De Luca and I. Torcicollo,}  \newblock Longtime behavior of vertical throughflows for binary mixtures in porous layers,  \emph{ Int. Journal of Non-Linear Mechanics} {\bf 52} (2013) 1--7.
\bibitem{CDDT} \newblock{F. Capone, V. De Cataldis,  R. De Luca and I. Torcicollo,}  \newblock{On the stability of vertical constant throughflows for binary mixtures in porous layers,} \emph{Int. Journal of Non-Linear Mechanics}  {\bf 59} (2014) 1--8.
%\bibitem{CGH} \newblock{F. Capone, M. Gentile and A.A. Hill,} \newblock{Penetrative convection in anisotropic porous media with    variable permeability,} \emph{Acta Mechanica} {\bf 216}  (2011), 49--58.
\bibitem{Cou}
% Researches into the Mathematical Principles of the Theory of Wealth, Hachette, Paris, France, 1838.
\newblock{A. Cournot,} \newblock{Recherches sur les Principes Mathematiques de la Theorie des Richesses} \newblock{Hachette, Paris. 1838}.
\bibitem{DA} \newblock{M. De Angelis, }  \newblock{Asymptotic Estimates Related to an integro Differential Equation, }  \emph{Nonlinear Dynamics and Systems Theory}  {\bf 13} (3)  (2013) 217--228.
%\bibitem{EAE} \newblock{E.M. Elabbasy, H. N. Agiza and A.A. Elsadany,}  \newblock{ Analysis of nonlinear triopoly  game with heterogeneous players,} \emph{Computers and Mathematics with Applications},  {\bf 57}  (2009), 488--499.
\bibitem{E} \newblock{A.A. Elsadany, }  \newblock{Dynamics of a delayed duopoly game with bounded rationality,}  \emph{ Mathematical and Computer Modelling} {\bf 52} (2010), 1497--1489.
\bibitem{Dir} \newblock{J. N., Flavin and S. Rionero,}  \newblock{Cross-diffusion influence on the nonlinear $l^2$-stability analysis for a Lotka-Volterra
reaction-diffusion system of PDEs,}  \emph{IMA  J. of Applied Mathematics}
{\bf 72} (2007) 540--555.
\bibitem{10} \newblock{J.N. Flavin and S. Rionero}, \emph{Qualitative estimates for partial differential equations: an introduction,}  \newblock{CRC Press,  Boca Raton (FL), 1996. }
\bibitem{Fri}  \newblock{R. Frish,}  \newblock{Monopoly- polipoly - the concept of force in the economy,} \newblock{International  Economic  Papers,} {\bf 1} ([1933] 1951) 23--26. 
\bibitem{G} \newblock{N. Giocoli, } \newblock{The Escape from conjectural variations: the consistency condition in duopoly theory from Bowley to Fellner,}  \emph{Cambridge Journal  of Economics,} {\bf 29 }(2005) 601--618.
\bibitem{K} \newblock{M. Kopel, } \newblock{Simple and complex adjustment dynamics in Cournot 
duopoly models,} \emph{ Chaos Solitons Fractals }  {\bf12}  (1996) 2031--2048.
%\bibitem{KTS} {P.K. Kythe, P. Puri and M.R. Schaferkotter,} \emph{PDEs in the Boundary Value Problems. }Second Edition. Chapman \& Hall/CRC,  2003.
%\bibitem{rjmaa} \newblock{G. Mulone and S. Rionero,} \newblock{On the nonlinear stability of the rotating B\'enard problem via the Lyapunov direct method,}  \emph{ Journal of Mathematical Analysis and Applications}  {\bf144} (1) (1989), 109--127.
%\bibitem{absence} \newblock{S. Rionero, }  \newblock{Absence of subcritical instabilities and global nonlinear stability for porous ternary diffusive-convective fluid mixtures,}  \emph{Phys. Fluids}, {\bf 24}  (2012), 104101.
\bibitem{rion} \newblock{S. Rionero,} \newblock{$L^2$-energy stability via new dependent variables for circumventing strongly nonlinear reaction terms,} \emph{Nonlinear Analysis }  {\bf 70}  (2009), 2530--2541.
\bibitem{rioneroult} \newblock{S. Rionero,} \newblock{$L^2$-energy decay of convective nonlinear PDEs reaction-diffusion systems via auxiliary ODEs systems,} \emph{Ricerche di Matematica,} DOI  10.1007/s11587-015-0231-2 (2015).
\bibitem{ternary}  \newblock{S. Rionero,}  \newblock{Stability of ternary reaction-diffusion dynamical systems,}  \emph{Rend. Lincei Mat. Appl.} {\bf 22} (2011) 245--268.
\bibitem{R} \newblock{S. Rionero,}  \newblock{A rigorous reduction of the $L^2$-stability  of the solutions to a nonlinear binary reaction-diffusion system of PDE's
to the stability of the solutions to a linear system of ODE's,}  \emph{Journal of 
Mathematical Analysis and Applications}, {\bf 319} (2006) 377--397.
%\bibitem{mcr} \newblock{S. Rionero,}  \newblock{Multicomponent diffusive-convective fluid motions in porous layers: Ultimately boundedness, absence of subcritical instabilities, and global nonlinear stability for any number of salts,}  \emph{Physics of Fluids} 2{\bf 5} (5) (2013), 1--25.
%\bibitem{heatmass} \newblock{S. Rionero,} \newblock{ Heat and mass transfer by convection in multicomponent Navier-Stokes mixtures: absence of subcritical instabilities and global nonlinear stability via the Auxiliary System Method,}  \emph{ Rend. Accademia Nazionale dei Lincei } {\bf 25} (2014), 1--44.
\bibitem{Ri} \newblock{S. Rionero,}  \newblock{A nonlinear $L^2$-stability  analysis for two species dynamics with dispersal, } \emph{ Math. Biosc. Eng.}  {\bf 3} (1)   (2006) 189--204.
%\bibitem{Rob} \newblock{S. Rionero,}  \newblock{Long time behavior of three competing species and mutualistic communities, } \emph{Asymptotic Methods in Nonlinear Wave Phenomena,}  \newblock{ World Scientific (2006) 171-186.}
 %\bibitem{Sor} \newblock{S. Rionero,} \newblock{Soret effects on the onset of convection in rotating porous layers via the \lq \lq auxiliary system method",}  \emph{Ricerche di Matematica} {\bf 62} (2) (2013), 183--208.
 \bibitem{RT1} \newblock{S. Rionero and I. Torcicollo,}  \newblock{Stability of a continuous reaction-diffusion Cournot-Kopel Duopoly Game Model,} \emph{ Acta Applicandae Mathematicae} {\bf 132}  (1), (2014) 505--513.  .
 \bibitem{RT} \newblock{S. Rionero and  I. Torcicollo,}  \newblock{On an ill-posed problem in nonlinear heat conduction, }
\emph{Transport Theory and Statistical Physics}  {\bf 29}  (1\&2)  (2000) 173--186.
%\bibitem{RT4} \newblock S. Rionero,  I. Torcicollo, \newblock On a nonlinear heat conduction problem, \newblock \emph{Proceedings of STAMM 2000, } (Editors: O'Donoghue, Flavin)  (2000), 178--184.
\bibitem{RT3} \newblock S. Rionero and I. Torcicollo, \newblock{On the pointwise continuous
dependence of
 an approximate solution of a nonlinear heat conduction ill-posed problem, } 
\emph{ Rend. Acc. Sci. Fis. Mat.,} Napoli,  {\bf LXVII} - Serie IV, (2000) 169--179.
\bibitem{T} \newblock{I. Torcicollo, } \newblock{On the dynamics of a non-linear Duopoly game model,}  \emph{Int. Journal of Non-Linear Mechanics} {\bf 57} (2013) 799--805.
\bibitem{T1}  \newblock{I. Torcicollo, } \newblock{Su alcuni problemi di diffusione non lineare}, \emph{Bollettino della Unione Matematica Italiana } A3 (3) (2000) 407-410.
\bibitem{TV1} \newblock{I. Torcicollo and M. Vitiello, }  \newblock{A note on the nonlinear pointwise stability for the equation $u_t=\Delta F(u)$ in the exterior of a sphere,} \emph{ Rend. Acc. Sc. Fis. Mat.  } Napoli,  {\bf LXX} (2003) 111--117.
\bibitem{TV3} \newblock I. Torcicollo, M. Vitiello, \newblock On the nonlinear diffusion in the exterior of a sphere,
\newblock  \emph{ Proceedings 11th International Conference on Wave and Stability in Continuous Media,} Porto Ercole, Italy.
Book Editor(s): Monaco, R; Bianchi, MP; Rionero, S;  
 % DOI: 10.1142/9789812777331_0069
  (2002) 563--568.
%\bibitem{TV2} \newblock{I. Torcicollo and M. Vitiello, }  \newblock{On the nonlinear diffusion in the exterior of a sphere.}  \emph{Rend. Acc. Sci. Fis. Mat., } Napoli, {\bf LXVIII } (2001), 139 -- 146.

%\bibitem{9} D.R. Merkin, Introduction to the theory of Stability, Springer texts in Appl. Math. V. 24, 1997.
 %\bibitem{Contr}  R. S. Cantrell, C. Cosner, Spatial ecology via reaction-diffusion equations. John Wiley and Sons, Ltd, %Chichester, UK, 2003. 
\bibitem{YY} \newblock{W. Yu and Y. Yu,} \newblock {A dynamic duopoly model with bounded rationality based on constant conjectural variation},  \emph{Economic Modelling},  {\bf 37} (1) (2014) 103--112. 
%\bibitem{ZZ} \newblock{L. Zhao and J. Zhang, } \newblock{Analysis of a duopoly game with heterogeneous players participating in carbon emission trading,}  \emph{Nonlinear Analysis: Modelling and Control} {\bf 19} (2014), 118--131.

 \end{thebibliography}
\end{document}